\documentclass[11pt]{article}
\pdfoutput=1
\usepackage{t1enc}
\usepackage[utf8]{inputenc}
\usepackage{amsmath,amssymb,amsfonts,amsthm,mathrsfs}
\DeclareFontFamily{OT1}{pzc}{}
\DeclareFontShape{OT1}{pzc}{m}{it}{<-> s * [1.2] pzcmi7t}{}
\DeclareMathAlphabet{\mathpzc}{OT1}{pzc}{m}{it}
\usepackage{graphicx}
\usepackage[all,ps]{xy}
\usepackage[colorlinks,urlcolor=cyan,citecolor=blue,linkcolor=blue]{hyperref}

\def\Embfrm{\mathop{\textstyle\rm Emb^{fr}_{\it m}}}
\def\max{\mathop{\textstyle\rm max}}
\def\d{\mathrm d}

\newenvironment{prf}%
{\par\noindent\textbf{Proof.\enspace\ignorespaces}}%
{~$\square$\par\medskip}%
\newenvironment{ack}%
{\medskip\par\noindent\textbf{Acknowledgement.\enspace\ignorespaces}}%
{\par}%
\newtheorem{thm}{Theorem}[section]%
\newtheorem{claim}[thm]{Claim}%
\newtheorem{crly}[thm]{Corollary}%
\theoremstyle{definition}
\newtheorem{defi}[thm]{Definition}%
\newtheorem{rmk}[thm]{Remark}%
\title{Stable Pontryagin--Thom construction for \\ proper maps}
\author{András Csépai}
\date{}

\begin{document}

\maketitle

\begin{abstract}
We will present proofs for two conjectures stated in \cite{rot}. The first one is that for an arbitrary manifold $W$, the homotopy classes of proper maps $W\times\mathbb{R}^n\to\mathbb{R}^{k+n}$ stabilise as $n\to\infty$, and the second one is that in a stable range there is a Pontryagin--Thom type bijection for proper maps $W\times\mathbb{R}^n\to\mathbb{R}^{k+n}$. The second one actually implies the first one and we shall prove the second one by giving an explicit construction.
\end{abstract}

\section{Introduction}

The main problem we are about to solve is to show the existence of a bijection between a set of homotopy classes and a set of cobordism classes. The first theorem of this type was proved by Pontryagin, who used it to calculate the homotopy groups of spheres. Pontryagin's idea was later generalised by Thom, but the general construction of the bijection remained almost the same. Unfortunately this Pontryagin--Thom construction does not entirely work in our setting, but we will see that a bit more complicated yet similar method does.

Let us first describe the objects we shall be working with.

\begin{defi}
A continuous map $f\colon X\to Y$ is said to be proper if $f^{-1}(C)$ is compact for all compact subsets $C\subset Y$. Two proper maps, $f,g\colon X\to Y$ are called proper homotopic, if there is a proper map $H\colon X\times[0,1]\to Y$ so that $H(\cdot,0)=f$ and $H(\cdot,1)=g$. The proper homotopy classes of proper maps $X\to Y$ will be denoted by $[X,Y]_{\mathrm{prop}}$.
\end{defi}

If $f\colon X\to Y$ is proper, then it is easy to see that the suspension of $f$ defined by
$$Sf\colon X\times\mathbb{R}\to Y\times\mathbb{R};~(x,t)\mapsto(f(x),t)$$
is also proper. Of course this construction can be defined for homotopies as well, so the suspensions of proper homotopic maps are also proper homotopic. Therefore there is a suspension map
$$S\colon[X,Y]_{\mathrm{prop}}\to[X\times\mathbb{R},Y\times\mathbb{R}]_{\mathrm{prop}}.$$

It is proved in \cite{rot} that for any vector bundle $E$ over a finite CW complex and for any sufficiently large $n\in\mathbb{N}$, the suspension
$$S\colon[E\times\mathbb{R}^n,\mathbb{R}^{k+n}]_{\mathrm{prop}}\to[E\times\mathbb{R}^{n+1},\mathbb{R}^{k+n+1}]_{\mathrm{prop}}$$
is a bijection, where $k\ge2$. We shall prove that the same is true if we take any (open or closed) manifold $W$ in the place of $E$ (and the $k\ge0$ condition is enough). This is the first conjecture in \cite{rot}.

Throughout this paper we will always mean smooth manifold when we say manifold and we will always assume that the manifold denoted by $W$ is connected.

\begin{defi}
Let $W$ be a manifold of dimension $m+k$ and $M\subset W$ an $m$-dimensional closed (compact) embedded submanifold with a trivial normal bundle. We say that a framing of $M$ is
$$\mathpzc{u}=(u_1,\ldots,u_k)\colon M\to NM$$
where $NM$ denotes the normal bundle of $M$ in $W$ and the $u_i$'s are pointwise linearly independent smooth normal vector fields. The pair $(M,\mathpzc{u})$ is called a framed submanifold.
\end{defi}

\begin{defi}
For an $(m+k)$-dimensional manifold $W$, a framed submanifold $(M,\mathpzc{u})$ is said to be framed cobordant to $(M',\mathpzc{u}')$ if there is a framed compact submanifold with boundary $(C,\mathpzc{v})$ of $W\times[0,1]$ so that $\partial C\subset W\times\{0,1\}$, $\partial C\cap W\times\{0\}=M$, $\partial C\cap W\times\{1\}=M'$ and the framings $\mathpzc{u}$ and $\mathpzc{u}'$ coincide with the restrictions of $\mathpzc{v}$. We will denote by $\Embfrm(W)$ the set of framed cobordism classes of $m$-dimensional closed manifolds embedded in $W$ with a framing.
\end{defi}

Our main result is the following:

\begin{thm}\label{t}
For any $m,k\in\mathbb{N}$ there is an $n_0\in\mathbb{N}$ so that for all manifolds $W$ of dimension $m+k$ and any $n>n_0$, there is a bijection
$$\Embfrm(W\times\mathbb{R}^n)\leftrightarrow[W\times\mathbb{R}^n,\mathbb{R}^{k+n}]_{\mathrm{prop}}.$$
\end{thm}

It is relatively easy to define a framed cobordism class for a given proper homotopy class; it is the same as in the standard Pontryagin construction, namely we take the preimage of a regular value. However, constructing the inverse of this map is a bit trickier; if we do not take $n\to\infty$, so if $n$ is not large enough, then it is not even always true, as it can be seen from the counterexamples given in \cite{rot}.

This bijection can be called stable Pontryagin--Thom construction and it is the second conjecture in \cite{rot}. We will prove the the first conjecture as a corollary of this theorem. It also implies that in a stable range $[W\times\mathbb{R}^n,\mathbb{R}^{k+n}]_\mathrm{prop}$ is in bijection with the 
%based 
homotopy classes of maps 
%sending the infinity to the infinity in 
between the one-point compactifications of the spaces, which we will denote by $[(W\times\mathbb{R}^n)^*,S^{k+n}]$.

Of course every proper map $X\to Y$ extends to a map $X^*\to Y^*$ by sending infinity to infinity, but the other way is not always true (we cannot get any map $X^*\to Y^*$ as the extension of a map $X\to Y$). By \cite{rot}, it is true for example for proper maps $E\to\mathbb{R}^k$ where $E$ is a vector bundle over a finite CW complex and $k$ is large enough. As we mentioned above, we will again prove the analogue of this for a trivial bundle over an arbitrary (open or closed) manifold $W$.

\begin{ack}
I would like to thank András Szűcs for his help.
\end{ack}

\section{The proof of theorem \ref{t}}

Fix an $n\in\mathbb{N}$, at first without any further conditions. For a proper map $f\colon W\times\mathbb{R}^n\to\mathbb{R}^{k+n}$ we may suppose (up to proper homotopy) that $f$ is smooth and $0$ is a regular value (by Sard's theorem). Then $M_f:=f^{-1}(0)$ is an $m$-dimensional submanifold of $W\times\mathbb{R}^n$ and it is compact, because $f$ is proper. If we fix a well-oriented basis in $T_0\mathbb{R}^{k+n}$ (i. e. the orientation of it is the same as that of the standard basis) and denote by $\mathpzc{u}_f$ the pullback of these vectors in $TW|_{M_f}$ given by the isomorphism $\d f_p\colon N_pM_f\to T_0\mathbb{R}^{k+n}$ for all $p\in M_f$, then $\mathpzc{u}_f$ is a framing of $M_f$. The framed submanifold $(M_f,\mathpzc{u}_f)$ is called the Pontryagin manifold of $f$.

Now we show that the framed cobordism class of $(M_f,\mathpzc{u}_f)$ does not depend on the choices we made.

\begin{claim}
If $f$ is proper homotopic to $g$ and $0$ is a regular value of $g$ as well, then $(M_f,\mathpzc{u}_f)$ is framed cobordant to $(M_g,\mathpzc{u}_g)$.
\end{claim}

\begin{prf}
If $H\colon W\times\mathbb{R}^n\times[0,1]\to\mathbb{R}^{k+n}$ is a proper homotopy, then we can assume that $H$ is smooth and $0$ is a regular value of $H$ too. We can also assume that for a small $\varepsilon>0$, $H$ fixes $f$ in $[0,\varepsilon]$ and $g$ in $[1-\varepsilon,1]$. Now $H$ has a Pontryagin manifold in $W\times\mathbb{R}^n\times[0,1]$, and it is easy to check that $(M_H,\mathpzc{u}_H)$ is a framed cobordism between $(M_f,\mathpzc{u}_f)$ and $(M_g,\mathpzc{u}_g)$.
\end{prf}\vspace{-.3em}

\begin{claim}
If $x\in\mathbb{R}^{k+n}$ is another regular value of $f$ with another fixed well-oriented basis in $T_x\mathbb{R}^{k+n}$, then the framed submanifold $f^{-1}(x)$ is framed cobordant to $(M_f,\mathpzc{u}_f)$.
\end{claim}

\begin{prf}
There is a compactly supported $\tau\colon\mathbb{R}^{k+n}\to\mathbb{R}^{k+n}$ diffeomorphism so that $\tau(x)=0$, $\d\tau_x$ carries the chosen basisvectors of $T_x\mathbb{R}^{k+n}$ to those of $T_0\mathbb{R}^{k+n}$ and $\tau$ is isotopic to the identity through a compactly supported isotopy. Then $\tau\circ f$ is proper homotopic to $f$, $0$ is a regular value of $\tau\circ f$, and $f^{-1}(x)=M_{\tau\circ f}$ and the given framing on $f^{-1}(x)$ is the same as $\mathpzc{u}_{\tau\circ f}$. Therefore by the previous claim, we have got a framed cobordism again.
\end{prf}

So we have constructed a well-defined map
$$[W\times\mathbb{R}^n,\mathbb{R}^{k+n}]_{\mathrm{prop}}\to\Embfrm(W\times\mathbb{R}^n).$$
What is left is to construct the inverse for it. In this part of the proof we will need $n$ to be a large number. Later it will be convenient to have $W\times\mathbb{R}^{n+1}$ instead of $W\times\mathbb{R}^n$, so in the remaining part of the proof we will use this.

Let $n_0:=\max\{1,m-k+2\}$ so that if $n\ge n_0$, the maps of $(m+1)$-dimensional manifolds into $(m+k+n+1)$-dimensionals can be approximated by embeddings (by Whitney's theorem). Fix an $n\ge n_0$ and a framed submanifold $(M,\mathpzc{u})$ of $W\times\mathbb{R}^{n+1}$. Our aim is to construct a Pontryagin--Thom collapse map $f\colon W\times\mathbb{R}^{n+1}\to\mathbb{R}^{k+n+1}$ so that $f$ is proper and the Pontryagin manifold of $f$ is $(M,\mathpzc{u})$.

Fix a Riemannian metric on $W$ and let $\mathpzc{u}=(u_1,\ldots,u_{k+n+1})$, where the $u_i$'s are orthonormal and orthogonal to $T_pM$ for all $p\in M$. We can make this assumption without loss of generality because for an arbitrary framing we may use the Gram--Schmidt process pointwise to construct a framing of this form, and the framed submanifold we get by this process will be framed cobordant to the initial one because this Gram--Schmidt process is a smooth deformation.

For the sake of simplicity we will call the last coordinate line of $\mathbb{R}^{n+1}$ in $T(W\times\mathbb{R}^{n+1})$ vertical.

\begin{claim}\label{c}
$(M,\mathpzc{u})$ is framed cobordant to a framed manifold $(M',\mathpzc{u}')$ where the last normal vector, $u'_{k+n+1}$ is vertical.
\end{claim}

\begin{prf}
By the compression theorem in \cite{rs}, $(M,u_{k+n+1})$ can be deformed by an isotopy to a submanifold $(M',u'_{k+n+1})$, where the normal vector $u'_{k+n+1}$ is vertical. This isotopy is an isotopy of $W\times\mathbb{R}^{n+1}$ therefore it is also a deformation of the other vector fields, this way we get the vector fields $u'_1,\ldots,u'_{k+n}$ on $M$. These vector fields are pointwise linearly independent and also independent of $TM'$, therefore $\mathpzc{u}':=(u'_1,\ldots,u'_{k+n+1})$ is a framing of $M'$.

If we denote the isotopy by $\varphi\colon W\times\mathbb{R}^{n+1}\times[0,1]\to W\times\mathbb{R}^{n+1}$ and $M_t:=\varphi(M\times\{t\})$ for all $t\in[0,1]$, then of course we may assume that for a small $\varepsilon>0$ we have $M_t=M$ if $t\in[0,\varepsilon]$ and $M_t=M'$ if $t\in[1-\varepsilon,1]$. Then the isotopy takes $\mathpzc{u}$ into a framing on each $M_t$, and so
$$\bigcup_{t\in[0,1]}M_t\times\{t\}\subset W\times\mathbb{R}^{n+1}\times[0,1]$$
is a framed cobordism between $(M,\mathpzc{u})$ and $(M',\mathpzc{u}')$.
\end{prf}

Hence we may assume that $u_{k+n+1}$ was initially vertical. Since it is a normal vector field, the projection of $M$ to $W\times\mathbb{R}^n$ is an immersion and because of the dimension condition made above we may also assume that it is an embedding.

\begin{claim}
$(M,\mathpzc{u})$ is framed cobordant to a framed manifold $(M',\mathpzc{u}')$ where $M'\subset W\times\mathbb{R}^n\times\{0\}$.
\end{claim}

\begin{prf}
Let $M'$ denote the image of $M$ under the projection to $W\times\mathbb{R}^n\times\{0\}$. For any $p\in M$ the tangent space at $p$ can be associated in a natural way with the tangent space at the projected image of $p$. Let the vectors $u'_i$ in the tangent spaces of the projected image be the ones associated to the vectors $u_i$ this way, and put $\mathpzc{u}':=(u'_1,\ldots,u'_{k+n+1})$.

 %The last real coordiante lines of $W\times\mathbb{R}^{n+1}$ are geodesics and we can parallel translate the vectors of $\mathpzc{u}$ along them, and this way we get a framing $\mathpzc{u}'$ of $M'$.

For all $t\in[0,1]$, the submanifold
$$M'_t:=\{(p_0,(1-t)p_1)\in W\times\mathbb{R}^{n+1}\mid p_0\in W\times\mathbb{R}^n,p_1\in\mathbb{R},(p_0,p_1)\in M\}$$
can be endowed with the framing we get similarly to how we defined the framing $\mathpzc{u}'$ for $M'$. Then
$$\bigcup_{t\in[0,1]}M'_t\times\{t\}\subset W\times\mathbb{R}^{n+1}\times[0,1]$$
is a framed cobordism between $(M,\mathpzc{u})$ and $(M',\mathpzc{u}')$.
\end{prf}

Hence we can also assume that $M$ was initially in $W\times\mathbb{R}^n\times\{0\}$, the last normal vector fied $u_{k+n+1}$ is vertical and the $u_i$'s are orthonormal and orthogonal to $M$ pointwise.

Our plan is now to construct a "nice" neighbourhood of $M$, then with the help of this neighbourhood define a map of $W\times\mathbb{R}^{n+1}$ to $\mathbb{R}^{k+n+1}$, and then prove that this map satisfies every condition we need.

Choose a small $\varepsilon>0$ so that the exponential map restricted to the open disc of radius $\varepsilon$ in $T_p(W\times\mathbb{R}^{n+1})$ is a diffeomorphism for all $p\in M$ ($M$ is compact, therefore there is such an $\varepsilon$).

For all $p\in M$ denote by $D_0(p)\subset W\times\mathbb{R}^{n+1}$ the $(k+n)$-dimensional open disc of radius $\varepsilon$ around $p$ so that $T_pD_0(p)$ is orthogonal to $u_{k+n+1}(p)$ and $T_pM$ (remember that the codimension of $M$ is $k+n+1$, so this disc is well-defined). Define
$$U_0:=\underset{p\in M}{\bigcup}D_0(p),$$
so $U_0$ is a tubular neighbourhood of $M$ in $W\times\mathbb{R}^n\times\{0\}$.

\begin{center}
\begin{figure}[h!]
\centering\includegraphics[scale=0.2]{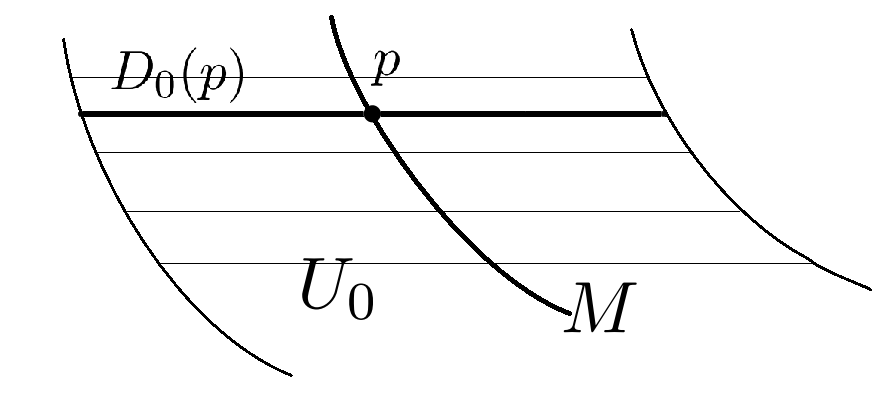}\label{1}
\begin{changemargin}{2cm}{2cm} 
\caption{\hangindent=1.4cm\small We indicated the case when $n=1$ and $k=0$. The arc in the middle represents $M$ and the strip around it is $U_0$.}
\end{changemargin} 
\end{figure}
\end{center}
\vspace{-.7cm}

For all $p\in M$ and $q\in D_0(p)$, put $q=(q_0,0)$ where $q_0\in W\times\mathbb{R}^n$. Let $t_q:=-\sqrt{\varepsilon^2-d(p,q)^2}$ (where $d(p,q)$ denotes their distance), so $t_q$ is the negative number for which $d(p,(q_0,t_q))=\varepsilon$. Using the notation $l(q):=\{q_0\}\times(t_q,\infty)$ define
$$U:=\bigcup_{p\in M}\bigcup_{q\in D_0(p)}l(q)=\bigcup_{p\in M}U_p,$$
where $U_p:=\underset{q\in D_0(p)}{\bigcup}l(q)$ is the fibre of $U$ above $p$.

We will also use the following notations: For all points $p\in M$ and $q\in D_0(p)$ we put $l_+(q):=\{q_0\}\times[0,\infty)$ and $l_-(q):=\{q_0\}\times(t_q,0]$. We define
$$D_+(p):=\bigcup_{q\in D_0(p)}l_+(q) ~\text{ and }~ D_-(p):=\bigcup_{q\in D_0(p)}l_-(q),$$
so $D_-(p)$ denotes the half of the $(k+n+1)$-dimensional open disc orthogonal to $M$ with centre $p$ and radius $\varepsilon$, in which the last coordinate of any point is non-positive.

\begin{center}
\begin{figure}[h!]
\centering\includegraphics[scale=0.2]{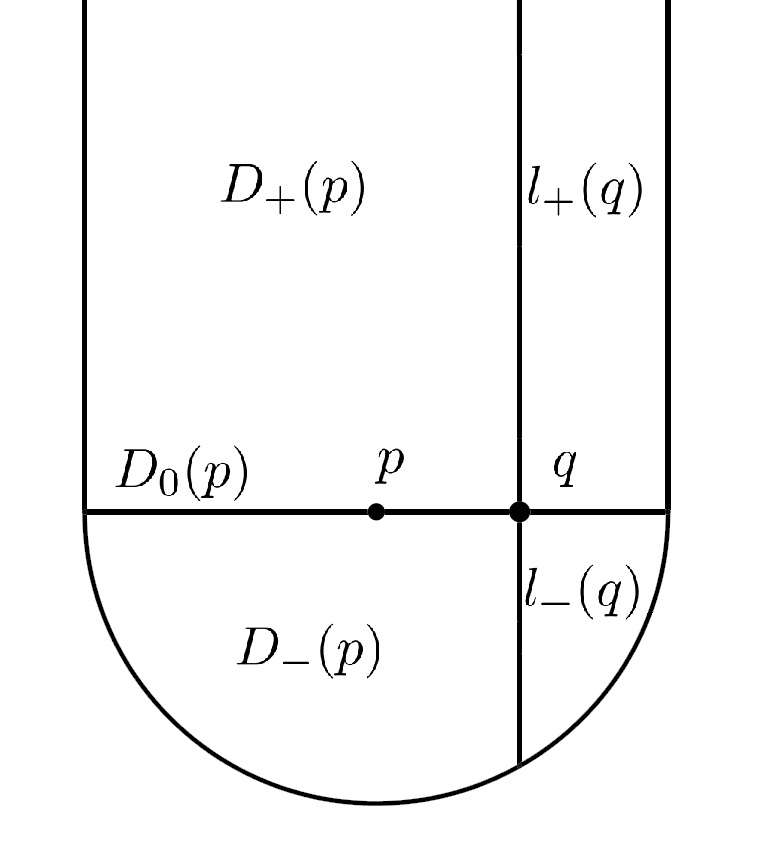}\label{2}
\begin{changemargin}{2cm}{2cm} 
\caption{\hangindent=1.4cm\small $U_p$ looks like this when $n=1$ and $k=0$. We also indicated $l(q)$ for a $q\in D_0(p)$.}
\end{changemargin} 
\end{figure}
\end{center}
\vspace{-.7cm}

\begin{rmk}
It is easy to see that $D_-(p)\cap D_+(p)=D_0(p)$ and $D_-(p)\cup D_+(p)=U_p$ for all $p\in M$ and $U_p$ is diffeomorphic to a neighbourhood of the ray ${0}\times[0,\infty)$ in $\mathbb{R}^{k+n+1}$. For different points $p$ the sets $U_p$ are disjoint because $\varepsilon$ was chosen sufficiently small, therefore their union, $U$ is a "nice" neighbourhood of $M$.
\end{rmk}

Now we are ready to define the desired Pontryagin--Thom collapse map for $M$. First fix a point $p\in M$ and define the map on $U_p$.

Let $v\in\mathbb{R}^{k+n+1}$ denote the last coordinate vector. Then there is a diffeomorphism
$$f_0\colon D_0(p)\to S^{k+n}\setminus\{-v\}$$
from the open disc $D_0(p)$ to the punctured sphere $S^{k+n}\setminus\{-v\}$ that maps the centre $p$ to the north pole $v$. We may of course choose $f_0$ so that the derivative $\d f_{0,p}$ maps the vectors $u_1(p),\ldots,u_{k+n}(p)$ to the standard basisvectors of the subspace $T_v(\mathbb{R}^{k+n}\times\{1\})\subset T_v\mathbb{R}^{k+n+1}$.

Then $f_0$ extends to a diffeomorphism
$$f_-\colon D_-(p)\to D^{k+n+1}\setminus\{-v\}.$$
This extension can be constructed in the following way: Take a diffeomorphism $\overline{D_-(p)}\approx D^{k+n+1}_-$, where $D^{k+n+1}_-:=\{(x_1,\ldots,x_{k+n+1})\in D^{k+n+1}\mid x_{k+n+1}\le0\}$; compose it with the quotient map
$$D^{k+n+1}_-\to D^{k+n+1}_-/(\partial D^{k+n+1}\cap D^{k+n+1}_-);$$
then associate the quotient space $D^{k+n+1}_-/(\partial D^{k+n+1}\cap D^{k+n+1}_-)$ with $D^{k+n+1}$ so that the image of the contracted boundary is $-v$. If we choose these diffeomorphisms so that the restriction of the composed map to $D_0(p)$ is $f_0$, then $f_-$ can be defined as the restriction of this map to $D_-(p)$.

\begin{center}
\begin{figure}[h!]
\centering\includegraphics[scale=0.2]{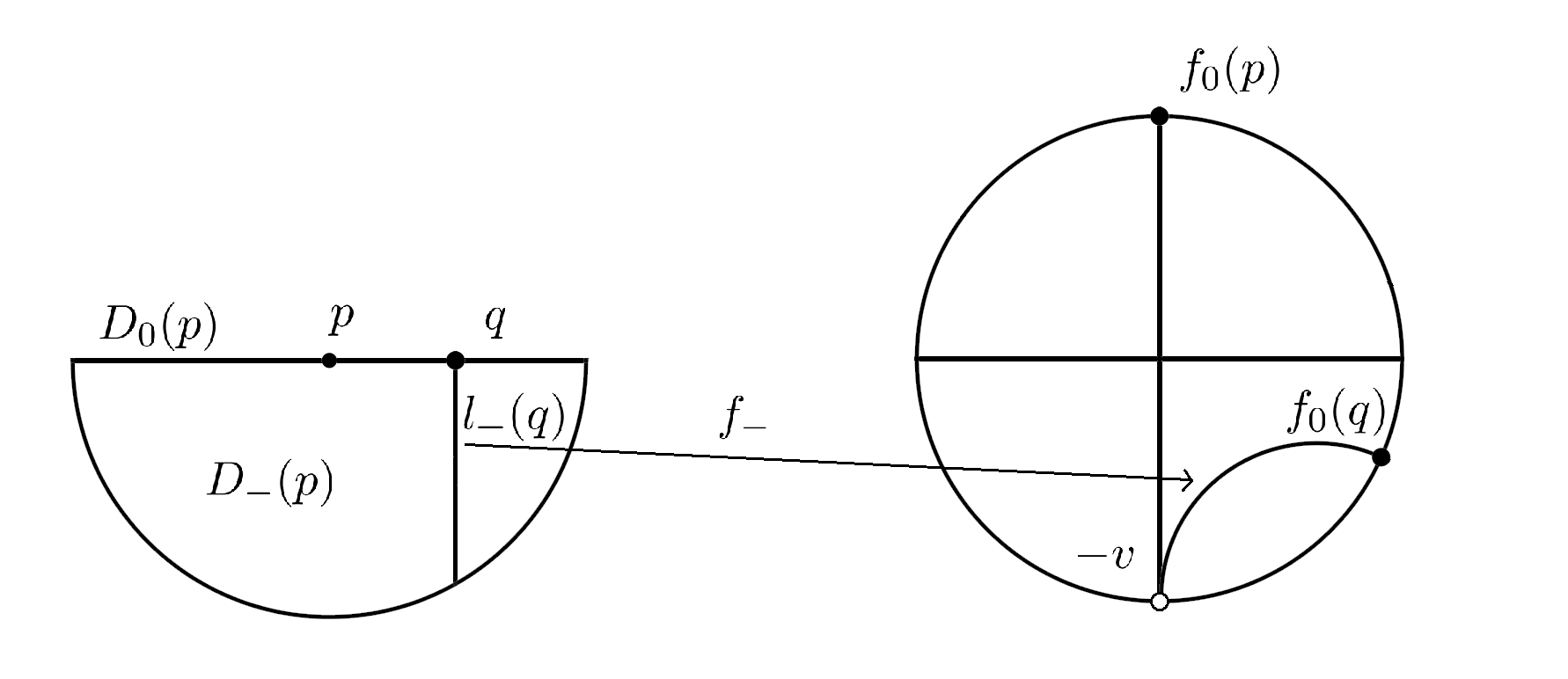}\label{3}
\begin{changemargin}{2cm}{2cm} 
\caption{\hangindent=1.4cm\small The left and the right sides of the figure represent $D_-(p)$ and $D^{k+n+1}\setminus\{-v\}$ respectively. The map $f_-$ restricted to $l_-(q)$ is also indicated.}
\end{changemargin} 
\end{figure}
\end{center}
\vspace{-.7cm}

In the following we will use the notation
$$[x\to):=\{(1+t)x\in\mathbb{R}^{k+n+1}\mid t\ge0\}$$
for all $x\in\mathbb{R}^{k+n+1}$, so $[x\to)$ is the ray from $x$ in the direction of $x$. 

For all $q\in D_0(p)$, the ray $l_+(q)$ has a linear bijection with $[f_0(q)\to)$ defined by $r\mapsto(1+d(r,p)-d(q,p))f_0(q)~(r\in l_+(q))$. The union of these maps is a diffeomorphism
$$f_+\colon D_+(p)\to\mathbb{R}^{k+n+1}\setminus((D^{k+n+1}\setminus S^{k+n})\cup[-v\to)).$$
Because of our conditions for $f_-$, the map $f_-\cup f_+$ is a diffeomorphism $U_p\approx\mathbb{R}^{k+n+1}\setminus[-v\to)$.

\begin{center}
\begin{figure}[h!]
\centering\includegraphics[scale=0.2]{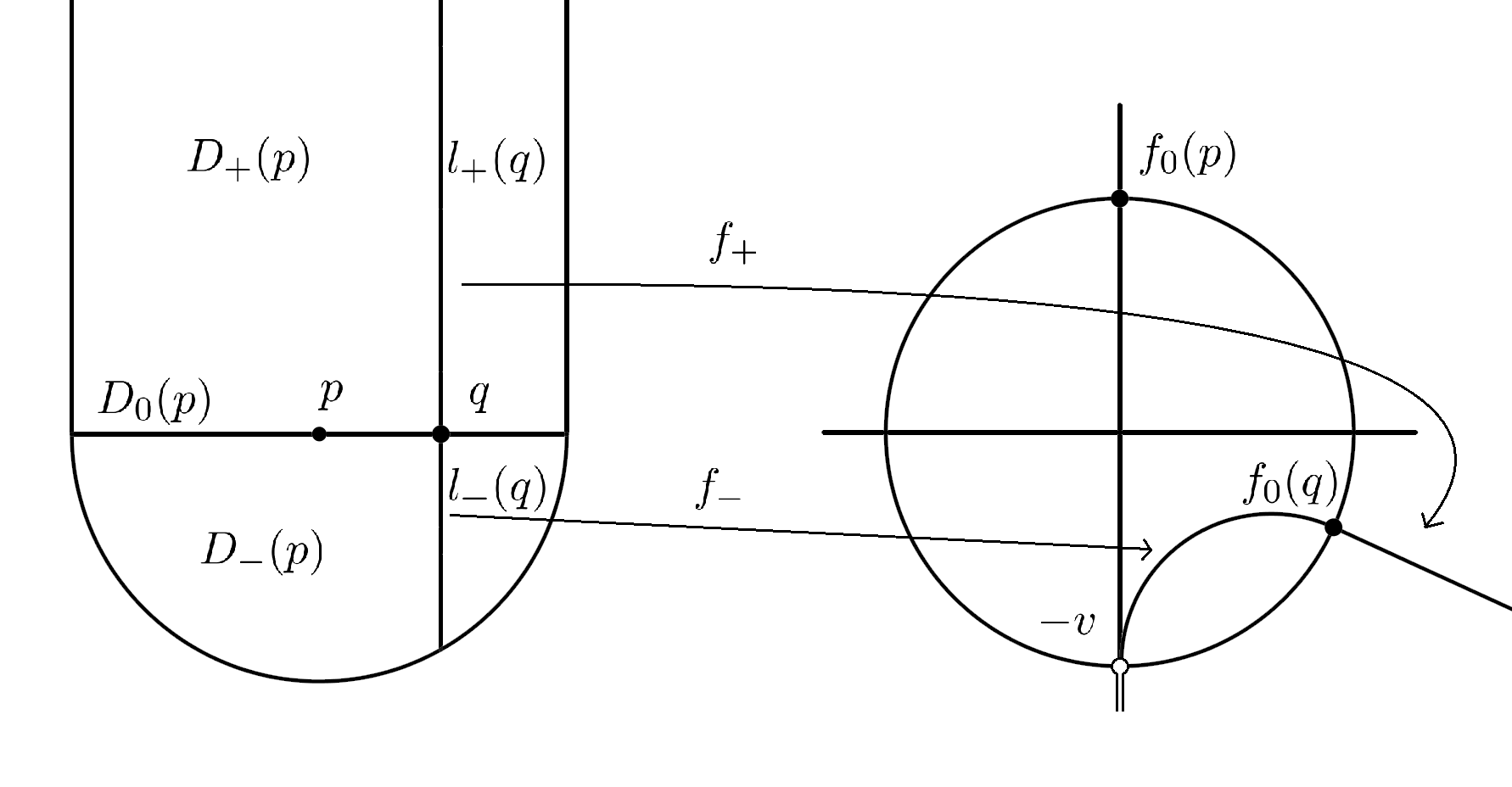}\label{a}
\begin{changemargin}{2cm}{2cm} 
\caption{\hangindent=1.4cm\small The left and the right sides of the figure represent $U_p$ and $\mathbb{R}^{k+n+1}\setminus[-v\to)$ respectively. The map $f_-\cup f_+$ restricted to $l(q)$ is also indicated.}
\end{changemargin} 
\end{figure}
\end{center}
\vspace{-.7cm}

The framing $\mathpzc{u}$ gives a natural diffeomorphism $\tau\colon U\to M\times U_p$ which maps each fibre $U_q$ into $\{q\}\times U_p$ isometrically. Denote the projection $M\times U_p\to U_p$ by $\pi$ and let
$$f_1:=(f_-\cup f_+)\circ\pi\circ\tau\colon U\to\mathbb{R}^{k+n+1}\setminus[-v\to).$$
Then the previous comment implies that the restriction of $f_1$ to an arbitrary fibre of $U$ is a diffeomorphism, so we only need to extend it to a map of $W\times\mathbb{R}^{n+1}$. For all $q\in W\times\mathbb{R}^{n+1}\setminus U$, the distance $d(q,M)$ from $M$ is well-defined and at least $\varepsilon$. Let
$$f_2\colon W\times\mathbb{R}^{n+1}\setminus U\to[-v\to);~q\mapsto-(1+d(q,M)-\varepsilon)v.$$

Define the Pontryagin--Thom collapse map as
$$f:=(f_1\cup f_2)-v\colon W\times\mathbb{R}^{n+1}\to\mathbb{R}^{k+n+1};~q\mapsto(f_1\cup f_2)(q)-v.$$

Now we have to prove that $f$ is a proper map for which $(M_f,\mathpzc{u}_f)=(M,\mathpzc{u})$. To see this, we first notice that $f$ restricted to any fibre $U_q$ of $U$ is the same as $f$ restricted to $U_p$ up to an isometric diffeomorphism. Therefore if we fixed $q$ at the beginning and did the same constructions using $q$ instead of $p$, then we would have constructed the same function $f$, so we can forget about the fixed point $p$.

\begin{claim}
$f$ is continuous.
\end{claim}

\begin{prf}
First we observe that for all $q\in U_0$ there is a unique $p\in M$ so that $q\in D_0(p)$ and for the point $p$ it holds that $d(q,M)=d(q,p)$. This is because $\varepsilon$ was chosen so that the exponential is a diffeomorphism on the $\varepsilon$ disc for all points of $M$ and $D_0(p)$ is orthogonal to $M$. Then the same is true if $q\in U_p$ is arbitrary, because if $q=(q_0,q_1)$, where $q_0\in W\times\mathbb{R}^n$ and $q_1\in\mathbb{R}$, then
$$d(q,M)=\sqrt{d(q_0,M)^2+q_1^2}=\sqrt{d(q_0,p)^2+q_1^2}=d(q,p).$$
This also implies that if $q\in\overline{U}_p$, then $d(q,M)=d(q,p)$ because the distance is continuous.

It is easy to see that $f_1$ and $f_2$ are both continuous, so we only need to prove that $f$ is continuous in the points of $\partial U$. Choose an arbitrary point $q\in\partial U$ and a sequence $(q_l)$ in $U$ so that $q_l\to q~(l\to\infty)$. We want to show that the sequence $(f(q_l))$ converges to $f(q)$, or equivalently $(f_1(q_l))$ converges to $f_2(q)$.

There is a $p\in M$ so that $q\in\partial U_p$. If $q\in\partial D_-(p)$, then $d(q,M)=d(q,p)=\varepsilon$, therefore $f_2(q)=-v$. Because of the construction of $f_-$ as a quotient map, $f_1(q_l)\to-v~(l\to\infty)$, hence $(f_1(q_l))$ indeed converges to $f_2(q)$.

If $q\in\partial D_+(p)$, then we may assume that all of the $q_l$'s are in $\underset{r\in M}{\bigcup}D_+(r)$. The sequence of the unit vectors $\frac{f_1(q_l)}{\lVert f_1(q_l)\rVert}$ converges to $-v$ because $(q_l)$ converges to $\partial U$. By the definition of $f_+$, the norm $\lVert f_1(q_l)\rVert=1+d(q_l,M)-d(q_{l,0},M)$, where $q_{l,0}$ denotes the component of $q_l$ in $W\times\mathbb{R}^n$. $(q_l)$ converges to $q\in\partial U$, therefore $d(q_{l,0},M)\to\varepsilon$ and $d(q_l,M)\to d(q,M)$ as $l\to\infty$. Hence the sequence of the norms $(\lVert f_1(q_l)\rVert)$ converges to $1+d(q,M)-\varepsilon$ and so $(f_1(q_l))$ converges again to $f_2(q)$.

Since we have proved this convergence for an arbitrary point and an arbitrary sequence, the continuity of $f$ follows.
\end{prf}

The map $f$ is smooth in a neighbourhood of $M$, $f^{-1}(0)=M$ and the framing we get as the pullback of the standard basisvectors is $\mathpzc{u}$. Therefore if we prove that $f$ is proper, then we get the desired result, that $(M_f,\mathpzc{u}_f)=(M,\mathpzc{u})$. %We only need to show that $f$ is proper.

\begin{claim}
$f$ is proper.
\end{claim}

\begin{prf}
Let $C\subset\mathbb{R}^{k+n+1}$ be an arbitrary compact subset. Then $C$ is closed, so $f^{-1}(C)$ is closed. $C$ is also bounded, therefore of course $C+v$ is bounded and $C+v=C_1\cup C_2$ where $C_1\subset\mathbb{R}^{k+n+1}\setminus[-v\to)$ and $C_2\subset[-v\to)$. The set $f_1^{-1}(C_1)$ is bounded, because $C_1$ is bounded and $f_1=(f_-\cup f_+)\circ\pi\circ\tau$ where $(f_-\cup f_+)^{-1}$ only increases the distance of two points by less than $2\varepsilon$, $\pi^{-1}$ is the product with a compact space (which is $M$) and $\tau^{-1}$ (the diffeomorphism between $M\times U_p$ and $U$) does not change the metric. $f_2^{-1}(C_2)$ is trivially bounded, because $C_2$ is bounded. Hence $f^{-1}(C)$ is a closed and bounded subset of $W\times\mathbb{R}^{n+1}$. If we assume that the Riemannian metric on $W$ is complete, then $f^{-1}(C)$ is compact by the Hopf--Rinow theorem. But we can assume that we have chosen the metric this way because of the results in \cite{no}.
\end{prf}

Now the only thing left to prove is that the Pontryagin--Thom construction we have defined is indeed well-defined. This means that the proper homotopy class of $f$ only depends on the framed cobordism class of $(M,\mathpzc{u})$.

\begin{claim}
If $(M,\mathpzc{u})$ is framed cobordant to $(M',\mathpzc{u}')$, then $f$ is proper homotopic to $f'$, the Pontryagin--Thom collapse map of $(M',\mathpzc{u}')$
\end{claim}

\begin{prf}
If $(C,\mathpzc{v})$ is a framed cobordism between $(M,\mathpzc{u})$ and $(M',\mathpzc{u}')$, then by the dimension condition for $n$, all of the constructions made to define $f$ for $(M,\mathpzc{u})$ have an analogue for $(C,\mathpzc{v})$. Therefore there is also a Pontryagin--Thom collapse map for $(C,\mathpzc{v})$, and it is easy to see that it is a proper homotopy between $f$ and $f'$.
\end{prf}

Hence we have an inverse map
$$\Embfrm(W\times\mathbb{R}^{n+1})\to[W\times\mathbb{R}^{n+1},\mathbb{R}^{k+n+1}]_{\mathrm{prop}}$$
for $n\ge n_0$, and our proof is complete.

\section{Corollaries}

Using theorem \ref{t} it is quite easy to solve the other problem in \cite{rot}, as we will now present.

\begin{crly}
For any $m,k\in\mathbb{N}$ there is an $n_0\in\mathbb{N}$ so that for all manifolds $W$ of dimension $m+k$ and any $n>n_0$, the suspension
$$S\colon[W\times\mathbb{R}^n,\mathbb{R}^{k+n}]_{\mathrm{prop}}\to[W\times\mathbb{R}^{n+1},\mathbb{R}^{k+n+1}]_{\mathrm{prop}}$$
is a bijection.
\end{crly}

\begin{prf}
The number $n_0$ will be the same as that of theorem \ref{t}. Then $[W\times\mathbb{R}^n,\mathbb{R}^{k+n}]_{\mathrm{prop}}$ is in bijection with $\Embfrm(W\times\mathbb{R}^n)$ and $[W\times\mathbb{R}^{n+1},\mathbb{R}^{k+n+1}]_{\mathrm{prop}}$ with $\Embfrm(W\times\mathbb{R}^{n+1})$. The suspension map for these is
\begin{alignat*}2
S\colon\Embfrm(W\times\mathbb{R}^n) & \to\Embfrm(W\times\mathbb{R}^{n+1})\\
[M,\mathpzc{u}] & \mapsto[M\times\{0\},(\mathpzc{u},u)],
\end{alignat*}
where $u$ denotes the constant vertical unit vector field and $[M,\mathpzc{u}]$ denotes the framed cobordism class of $(M,\mathpzc{u})$.

Claim \ref{c} implies that every framed cobordism class in $\Embfrm(W\times\mathbb{R}^{n+1})$ has a representative of the form $(M',(\mathpzc{u}',u))$. The projection of $M'$ to $W\times\mathbb{R}^n$ is an immersion and because of the dimension condition we may also assume that it is an embedding. If we denote its image under the projection by $M$ and the projected framing by $\mathpzc{u}$, then $(M,\mathpzc{u})$ is a framed submanifold of $W\times\mathbb{R}^n$ for which $(M',(\mathpzc{u}',u))$ is framed cobordant to $(M\times\{0\},(\mathpzc{u},u))$. Hence $S$ is surjective.

But claim \ref{c} can also be used to project framed cobordisms in $W\times\mathbb{R}^{n+1}\times[0,1]$ into $W\times\mathbb{R}^n\times[0,1]$, thus $S$ is also injective.
\end{prf}

Another nice corollary of theorem \ref{t} is that the proper homotopy classes can completely be classified by the homotopy classes between the one-point compactifications.

\begin{crly}
For any $m,k\in\mathbb{N}$ there is an $n_0\in\mathbb{N}$ so that for all manifolds $W$ of dimension $m+k$ and any $n>n_0$, there is a bijection
$$[W\times\mathbb{R}^n,\mathbb{R}^{k+n}]_{\mathrm{prop}}\leftrightarrow[(W\times\mathbb{R}^n)^*,S^{k+n}].$$
\end{crly}

\begin{prf}
If $n_0$ is again the same as before, then by theorem \ref{t} $[W\times\mathbb{R}^n,\mathbb{R}^{k+n}]_{\mathrm{prop}}$ is in bijection with $\Embfrm(W\times\mathbb{R}^n)$. And $[(W\times\mathbb{R}^n)^*,S^{k+n}]$ is in bijection with $[((W\times\mathbb{R}^n)^*,\infty),(S^{k+n},\infty)]$ (the based homotopy classes of maps sending $\infty$ to $\infty$), because for an arbitrary space $X$ an arbitrary map $f\colon X^*\to S^{k+n}$ can be combined with a rotation of the sphere that takes $f(\infty)$ into $\infty$ and this combined map will of course be homotopic to $f$. This is well-defined because $X^*\times[0,1]\cong(X\times[0,1])^*$, so the same works for homotopies.

Therefore it is enough to give a bijection between $\Embfrm(W\times\mathbb{R}^n)$ and $[((W\times\mathbb{R}^n)^*,\infty),(S^{k+n},\infty)]$. This is a standard Pontryagin--Thom type bijection and we will just sketch the proof which is a simpler version of the proof of theorem \ref{t}.

For a map $f\colon((W\times\mathbb{R}^n)^*,\infty)\to(S^{k+n},\infty)$ we can define the Pontryagin manifold in the same way as we did before. This is a framed closed submanifold which does not contain $\infty$, so it is in $W\times\mathbb{R}^n$. By the first part of the proof of theorem \ref{t}, this is a well-defined map
$$[((W\times\mathbb{R}^n)^*,\infty),(S^{k+n},\infty)]\to\Embfrm(W\times\mathbb{R}^n).$$

For a framed closed submanifold $(M,\mathpzc{u})$ of $W\times\mathbb{R}^n$, the framing gives a diffeomorphism between $U$, a tubular neighbourhood of $M$, and $M\times\mathbb{R}^{k+n}$. If we project this to $\mathbb{R}^{k+n}$, then we have a map $U\to\mathbb{R}^{k+n}$ and we can extend this to $(W\times\mathbb{R}^n)^*$ so that we map everything else to $\infty$. Similarly to the end of the proof of theorem \ref{t}, this is a well-defined inverse map
$$\Embfrm(W\times\mathbb{R}^n)\to[((W\times\mathbb{R}^n)^*,\infty),(S^{k+n},\infty)].$$

Combining these three bijections, we get the desired result.
\end{prf}


\begin{thebibliography}{00}

\bibitem{no} K. Nomizu, H. Ozeki, {\it The existence of complete Riemannian metrics}, Proc. Amer. Math. Soc. 12 (1961), 889-891
\bibitem{rot} T. O. Rot, {\it Homotopy classes of proper maps out of vector bundles}, preprint arXiv:1808.08073
\bibitem{rs} C. Rourke, B. Sanderson, {\it The compression theorem I}, Geometry and Topology 5 (2001), 399–429

\end{thebibliography}
\end{document}